\documentclass[11pt]{article}
\usepackage{amsmath}
\usepackage{dsfont}
\usepackage{mathrsfs}
\usepackage{amsmath,amssymb}
\usepackage{amsfonts}
\usepackage{hyperref}
\usepackage{amsthm}
\usepackage{graphicx}
\usepackage{subfigure}
\usepackage{xcolor}

\renewcommand{\qed}{\hfill\small{$\square$}\normalsize}

\hfuzz=\maxdimen
\theoremstyle{definition}
\newtheorem{lemma}{Lemma}[section]
\newtheorem{definition}[lemma]{Definition}
\newtheorem{proposition}[lemma]{Proposition}
\newtheorem{theorem}[lemma]{Theorem}
\newtheorem{corollary}[lemma]{Corollary}

\numberwithin{equation}{section}
\renewcommand{\proof}{\textbf{Proof. }}
\renewcommand{\qed}{\hfill\small{$\square$}\normalsize}

\DeclareFixedFont{\Acknowledgment}{OT1}{cmr}{bx}{n}{14pt}
\begin{document}

\title{\bf A Discrete Ricci Flow on Surfaces in Hyperbolic Background Geometry}
\author{Huabin Ge, Xu Xu}
\maketitle

\begin{abstract}
In this paper, we generalize our results in \cite{GX3} to triangulated surfaces in hyperbolic background geometry,
which means that all triangles can be embedded in the standard hyperbolic space.
We introduce a new discrete Gaussian curvature by dividing the classical discrete Gauss
curvature by an area element, which could be
taken as the
area of the hyperbolic disk packed at each vertex. We prove that the corresponding discrete Ricci flow converges if and only if there exists a circle packing metric with zero
curvature. We also prove that the flow converges if the initial curvatures are all negative.
Note that, this result does not require the existence of zero curvature metric or Thurston's combinatorial-topological condition.
We further generalize the definition of combinatorial curvature to any given area element and
prove the equivalence between the existence of zero curvature metric and the convergence of the corresponding
flow.
\end{abstract}

\section{Introduction}
Consider a compact surface $X$ with a triangulation $T=\{V,E,F\}$, where the symbols $V,E,F$ represent the set of vertices, edges and faces respectively. A positive function $r:V\rightarrow (0,+\infty)$ defined on the vertices  is called a circle packing metric and a function $\Phi: E\rightarrow [0, \pi/2]$  is called a weight on the triangulation. Throughout this paper, a function defined on vertices is regarded as a column vector and $N=V^{\#}$ is used to denote the number of vertices. Moreover, all vertices, marked by $v_{1},...,v_{N}$, are supposed to be ordered one by one and we often write $i$ instead of $v_i$. Thus we may think of circle packing metrics as points in $\mathds{R}^N_{>0}$, $N$ times of Cartesian product of $(0,\infty)$.

Given $(X,T,\Phi)$, every circle packing metric $r$ determines a piecewise linear metric on $X$ by attaching edge $e_{ij}$ a length $l_{ij}=\sqrt{r^2_i+r^2_j+2r_ir_jcos(\Phi_{ij})}$. This length structure makes each triangle in $T$ isometric to a Euclidean triangle. Furthermore, the triangulated surface $(X,T)$ is composed of many Euclidean triangles glued coherently. This case is called Euclidean background geometry in \cite{CL1}. If the length $l_{ij}$ of the edge is determined by the hyperbolic cosine law $\cosh l_{ij}=\cosh r_i\cosh r_j+\sinh r_i\sinh r_j\cos \Phi_{ij}$, then each triangle in $T$ is isometric to a hyperbolic triangle embedded in the standard hyperbolic space $\mathds{H}^2$. In this case the triangulated surface $(X,T)$ is composed of many hyperbolic triangles glued coherently, which is called hyperbolic background geometry.

We use $(X,T,\Phi,\mathds{H}^2)$ to denote the space we want to study in the following,
where $X$ is a closed surface, $T$ is a fixed triangulation of $X$, $\Phi$ is a fixed weight function defined on edges, and $\mathds{H}^2$ represents the hyperbolic background geometry. It was proved by Thurston \cite{T1} that, whenever $\{i,j,k\}\in F$, these three positive numbers $l_{ij}, l_{ik}, l_{jk}$ satisfy the triangle inequalities. Thus the combinatorial triangle $\{i,j,k\}$ with lengths $l_{ij}, l_{ik}, l_{jk}$ forms a hyperbolic triangle in $\mathds{H}^2$.

For the hyperbolic triangle $\triangle v_iv_jv_k$, the inner angle of this triangle at $v_i$ is
denoted by $\theta_i^{jk}$, and the classical combinatorial Gauss curvature $K_i$ at $v_i$ is defined to be
\begin{equation}\label{classical Gauss curvature}
K_i=2\pi-\sum_{\{i,j,k\}\in F}\theta_i^{jk}.
\end{equation}
Notice that $\theta_i^{jk}$ can be calculated by hyperbolic cosine law, thus $\theta_i^{jk}$ and $K_i$ are elementary functions of the circle packing metric $r$.

Obviously, discrete Gaussian curvatures $K$ is determined by circle packing metric $r$. Conversely, Thurston proved \cite{T1} that $K$-curvature map $r\rightarrow K$ is injective, which implies that the metric $r$ is determined by its $K$-curvature. Among all circle packing metrics, a metric with all curvatures zero, which is unique if it exists, is of special interest. Zero curvature metric does not always exist.
At first glance, the existence of zero curvature metric depends on the topological information of the surface. In fact, for the hyperbolic triangle $\triangle v_iv_jv_k$, we have $\theta_i^{jk}+\theta_j^{ik}+\theta_k^{ij}=\pi-Area(\triangle v_iv_jv_k)$. Using this formula, one get the combinatorial Gauss-Bonnet formula \cite {CL1}
\begin{equation}\label{Gauss-Bonnet formula}
\sum_{i=1}^NK_i=2\pi \chi(X)+Area(X),
\end{equation}
which implies, if there exists a circle packing metric with zero curvature, then $\chi(M)$ must be negative. However, only topological constraint is not enough. In his work on constructing hyperbolic metrics on $3$-manifolds, Thurston \cite{T1} discovered additional combinatorial-topological conditions which is equivalent to the existence of zero curvature metric.

\begin{proposition}\label{Lemma-Thurston's condition}
(\cite{T1}) Given $(X,T,\Phi,\mathds{H}^2)$, there exists a zero curvature metric if and only if the following two combinatorial-topological conditions are satisfied simultaneously:
\begin{enumerate}
  \item For any three edges $e_1, e_2, e_3$ forming a null homotopic loop in $X$, if $\sum_{i=1}^3\Phi(e_i)\geq\pi$, then $e_1, e_2, e_3$ form the boundary of a triangle of $T$;
  \item For any four edges $e_1, e_2, e_3, e_4$ forming a null homotopic loop in $X$, if $\sum_{i=1}^4\Phi(e_i)\geq2\pi$, then $e_1, e_2, e_3, e_4$ form the boundary of the union of two adjacent triangles.\\
\end{enumerate}
\end{proposition}

Bennett Chow and Feng Luo introduced the combinatorial curvature flow to study the problem. They defined a combinatorial Ricci flow
\begin{equation}\label{ChowLuo's flow}
\frac{dr_i}{dt}=-K_i\sinh r_i
\end{equation}
in \cite{CL1} and  proved that the flow (\ref{ChowLuo's flow}) converges if and only if there exists a metric with all curvatures are zero, if and
 only if Thurston's combinatorial-topological conditions are satisfied.

Inspired by \cite{CL1} and \cite{Ge1}, we studied the negative gradient flow of discrete Calabi energy,
\begin{equation}\label{GeXu's-hypbolic-Calabi-flow}
\frac{dr_i}{dt}=-\Delta K_i\sinh r_i,
\end{equation}
which is called combinatorial Calabi flow in \cite{Ge1}. We proved that combinatorial Calabi flow converges if and only if
zero curvature metric exists, assuming a uniform upper bound of metric $r$.

The paper is organized as follows. In Section \ref{Definitions of discrete Gaussian curvature and Ricci flow}, we introduce
a new definition of discrete Gauss curvature for triangulated surfaces in hyperbolic background geometry and
the corresponding combinatorial Ricci flow. In Section \ref{Convergence of Ricci flow when zero curvature metric exists},
we prove that the existence of zero curvature metric is equivalent to the convergence of combinatorial Ricci flow
we introduce in the paper, which provides another description of Thurston's existence of circle packing.
In Section \ref{Convergence with nonpositive initial curvature metrics}, we prove that the existence of metrics with nonpositive
curvatures ensures the existence of zero curvature metric.
In Section \ref{Combinatorial Ricci flow for general area element}, we generalized the definition of combinatorial
curvature to general area element and then prove that the existence of zero curvature metric in this sense is also
equivalent to the convergence of the corresponding combinatorial curvature flow.

\section{Definitions of discrete Gaussian curvature and Ricci flow}\label{Definitions of discrete Gaussian curvature and Ricci flow}
By the analysis in \cite{GX3}, the classical combinatorial Gaussian curvature is not a good candidate for approximating the classical Gaussian curvature. Following the idea and the approach in \cite{GX3}, we have the following definition of combinatorial Gaussian curvature.
\begin{definition}
Givin $(X,T,\Phi,\mathds{H}^2)$ with circle packing metric $r$, the modified combinatorial Gaussian curvature at the vertex $i$ is defined as
\begin{equation}\label{definition of modified Gauss curvature}
R_i=\frac{K_i}{4\pi\sinh^2\frac{r_i}{2}},
\end{equation}
where $K_i$ is the classical combinatorial Gaussian curvature given by (\ref{classical Gauss curvature}).
\end{definition}

In fact, the denominator $4\pi\sinh^2\frac{r_i}{2}$ is just the area of the hyperbolic disk with radius $r$ packed at vertex $i$.
To study the constant combinatorial curvature problem, we introduce a combinatorial Ricci flow which is an
analogue of the smooth Ricci flow.
\begin{definition}
Given $(X, T, \Phi,\mathds{H}^2)$ with circle packing metric $r$, the combinatorial Ricci flow is defined as
\begin{equation}\label{def-ricci flow}
\frac{dg_i}{dt}=-R_ig_i,
\end{equation}
where $g_i=\sinh^2\frac{r_i}{2}$.
\end{definition}

For the following application, we usually write the combinatorial Ricci flow (\ref{def-ricci flow}) in the following form
\begin{equation}\label{def-ricci flow2}
\frac{du_i}{dt}=-\widetilde{R}_i,
\end{equation}
where $\widetilde{R}_i=\frac{K_i}{2\pi\sinh^2r_i}$ is a modification of the curvature $R_i$
and $u_i=\ln \tanh\frac{r_i}{2}$ is a coordinate transformation,
which maps $\mathds{R}^N_{>0}$ homeomorphically onto $\mathds{R}^N_{<0}$.

The flow (\ref{def-ricci flow}) is in fact an ODE system.
Therefore, the solution always exists locally around the initial time $t=0$
since all coefficients are smooth and then locally Lipschitz continuous. It's easy to get
\begin{proposition}\label{prop-converg-imply-zccpm-exist}
If the hyperbolic combinatorial Ricci flow (\ref{def-ricci flow}) converges, then there exists a circle packing metric with $0$ $R$-curvature.
\end{proposition}
\proof
Suppose the solution $r(t)\rightarrow r^*$ as $t\rightarrow +\infty$. Then we have $g(t)\rightarrow g^*$, which implies that
there exists a sequence $\xi_n\rightarrow +\infty$ such that
$$g_i(n+1)-g_i(n)=g'_i(\xi_n)=-R_ig_i(\xi_n)\rightarrow 0$$
as $n\rightarrow +\infty$. As $r(t)\rightarrow r^*$, we have  $R(t)\rightarrow R^*$. Thus $R^*g^*=0$, which
implies $R^*=0$ and $r^*$ is a metric with zero curvature. \qed\\

Obviously, the $K$-curvature, the $R$-curvature and the $\widetilde{R}$-curvature are different.
However, a metric has zero $K$-curvature if and only if it has zero $R$-curvature, if and only if it has zero $\widetilde{R}$-curvature.
Thus in the following we needn't distinguish them and just call them ``\textbf{zero curvature metric}" for short.

To study the long time behaviors of combinatorial Ricci flow (\ref{def-ricci flow}), we shall use two different ways in the following subsections.

\section{Convergence of Ricci flow when zero curvature metric exists}\label{Convergence of Ricci flow when zero curvature metric exists}
In this subsection, we shall study carefully the relationship
between the existence of zero curvature metric and the long time behavior of the combinatorial Ricci flow.
We also derive some convergence results under the assumption of existence of a zero curvature metric,
which is necessary for the convergence of flow (\ref{def-ricci flow}) by Proposition \ref{prop-converg-imply-zccpm-exist}.
The main tools used here are standard in ODE theory.

The following two Lemmas are useful.
\begin{lemma}\label{Lemma-L-positive-definite}
(\cite{CL1} \cite{Gx1}) Given $(X, T, \Phi, \mathds{H}^2)$ with circle packing metric $r$. $u$ is the coordinate transformation of $r$ with $u_i=\ln \tanh\frac{r_i}{2}$. Then
$$L=\frac{\partial (K_1,\cdots,K_N)}{\partial(u_1,\cdots, u_N)}=A+L_B$$
is positive definite, where $A$ is a positive diagonal matrix and $L_B$ is positive semi-definite with rank $N-1$ and kernel $t(1,\cdots, 1)^T$ and
could be written as
\begin{gather*}
\left(L_B\right)_{ij}=
\begin{cases}
\,\,\sum\limits_{k \sim i}B_{ik} \,, & \text{$ j=i,$} \\
\,\,\,\,\,-B_{ij} \,,& \text{$ j\sim i,$} \\
\,\,\,\,\,\,\,\,\,0 \,,& \text{$ j\nsim i,\, j\neq i,$}
\end{cases}
\end{gather*}
with $B_{ij}>0$.
 \qed\\
\end{lemma}

\begin{lemma}\label{Lemma-edge-big-then-angle-zero}
(Lemma 3.5, \cite{CL1}) For any $\epsilon>0$, there exists a number $l$ so that when $r_i>l$, the inner angle $\theta_i$ in the hyperbolic triangle $\triangle v_iv_jv_k$ is smaller than $\epsilon$.

\proof
This result has been stated in \cite{CL1}, here we give a proof just for completeness.
We will prove that $\theta_{i}^{jk}\rightarrow 0$ uniformly as $r_i\rightarrow +\infty$.
By the hyperbolic cosine law, we have
\begin{equation}
\begin{aligned}
\cos \theta_i^{jk}
=&\frac{\cosh l_{ij}\cosh l_{ik}-\cosh l_{jk}}{\sinh l_{ij}\sinh l_{ik}}\\
=&\frac{\cosh(l_{ij}+l_{ik})+\cosh(l_{ij}-l_{ik})-2\cosh l_{jk}}{\cosh(l_{ij}+l_{ik})-\cosh(l_{ij}-l_{ik})}\\
=&\frac{1+\lambda-2\mu}{1-\lambda},
\end{aligned}
\end{equation}
where $\lambda=\frac{\cosh(l_{ij}-l_{ik})}{\cosh(l_{ij}+l_{ik})}$ and $\mu=\frac{\cosh l_{jk}}{\cosh(l_{ij}+l_{ik})}$.
Note that
$$0<\lambda<\frac{\cosh l_{ij}}{\cosh (l_{ij}+l_{ik})}<\frac{1}{\cosh l_{ik}}<\frac{1}{\cosh r_i},$$
we have $\lambda\rightarrow 0$ uniformly as $r_i\rightarrow +\infty$.
So we just need to prove that $\mu\rightarrow 0$ uniformly as $r_i\rightarrow +\infty$.
Note that
$$
\cosh l_{jk}
\leq \cosh r_j\cosh r_k+\sinh r_j\sinh r_k
=\cosh(r_j+r_k),
$$
we just need to prove that $l_{ij}+l_{ik}-(r_j+r_k)\rightarrow +\infty$ uniformly as $r_i\rightarrow +\infty$.
Note that
$$
\cosh l_{ij}
=\cosh r_i\cosh r_j+\sinh r_i\sinh r_j\cos \Phi_{ij}
\geq\cosh r_i\cosh r_j\geq \frac{1}{4}e^{r_i+r_j},
$$
we have
$$l_{ij}\geq \ln \cosh l_{ij}\geq r_i+r_j-\ln 4.$$
Similarly, we have
$$l_{ik}\geq \ln \cosh l_{ik}\geq r_i+r_k-\ln 4.$$
Then we have
$$l_{ij}+l_{ik}-(r_j+r_k)\geq 2r_i-2\ln 4\rightarrow +\infty$$
uniformly as $r_i\rightarrow +\infty$.
\qed\\
\end{lemma}

\begin{theorem}\label{Thm-0curv-metric-exist-imply converg}
Assuming there exists a zero curvature metric $r^*$ on $(X, T, \Phi,\mathds{H}^2)$. Then the solution $r(t)$ to the combinatorial Ricci flow (\ref{def-ricci flow}) exists for all time $t\in [0,+\infty)$ and converges exponentially fast to $r^*$.
\end{theorem}
\proof
Denote $u^*$ as the corresponding $u$-coordinate of $r^*$. Consider the combinatorial Ricci potential (first introduced by Chow and Luo in \cite{CL1}),
\begin{equation}\label{def-G-ricci-potential}
F(u)\triangleq\int_{u^*}^u\sum_{i=1}^N K_idu_i.
\end{equation}
Then $Hess(F)=L$ is positive definite. We had proved that $F(u)\geq F(u^*)=0$ and $u^*$ is the unique minimum point and zero point of $F$. Moreover, $\lim\limits_{\|u\|\rightarrow +\infty}F(u)=+\infty$ (Lemma B.1, \cite{Gx1}). Then
\begin{equation}
\frac{dF(u(t))}{dt}=\nabla F\cdot \frac{du}{dt}=\sum_iK_i\cdot (-\widetilde{R}_i)=-\sum_i\frac{K_i^2}{2\pi\sinh^2r_i}\leq 0.
\end{equation}
So $F$ is decreasing along the flow. Using $\lim\limits_{\|u\|\rightarrow +\infty}F(u)=+\infty$, we know $\{u(t)\}$ is bounded in $\mathds{R}^N_{<0}$, i.e. there exists $c>0$, such that $u_i(t)\geq -c$ for any $i\in V$ and $t\in [0,+\infty)$.
Therefore $r_i(t)\geq \delta >0$ for any $i\in V$ and $t\in [0,+\infty)$, where $\delta=\ln \frac{1+e^{-c}}{1-e^{-c}}$.

The combinatorial Ricci flow (\ref{def-ricci flow}) can be written as $(\cosh r_i)'(t)=-\frac{K_i}{2\pi}$. Note that $K_i$ is uniformly bounded by $(2-d)\pi<K_i<2\pi$, where $d$ is the maximal degree at vertices. Hence $\cosh r_i(t)\leq \cosh r_i(0)+\frac{d-2}{2} t$ and then $r_i(t)\leq \ln2+\ln\left(\cosh r_i(0)+\frac{d-2}{2}t\right)$ for each $i\in V$, which implies that the solution $r(t)$ exists for all time $t\in[0,+\infty)$.

We claim that $r(t)$ is uniformly bounded above.
If not, then there exists at least one vertex $i\in V$, such that $\overline{\lim\limits_{t\rightarrow +\infty}}r_i(t)=+\infty$.
For vertex $i$, using Lemma \ref{Lemma-edge-big-then-angle-zero}, we can choose $l>0$ large enough so that, whenever $r_i>l$,
the inner angle $\theta_i^{jk}$ is smaller than $\frac{\pi}{d_i}$, where $d_i$ is the degree at vertex $i$. Thus $K_i=2\pi-\sum \theta_i^{jk}>\pi$. Choose a time $t_0$ such that $r_i(t_0)>l$, this can be done since $\overline{\lim\limits_{t\rightarrow +\infty}}r_i(t)=+\infty$. Denote $a=\inf\{\,t<t_0\,|\,r_i(t)>l\,\}$, then $r_i(a)=l$. Let's look at what happens to flow (\ref{def-ricci flow}) in the interval $[a, t_0]$. Note that flow (\ref{def-ricci flow}) can be written as $r'_i(t)=-\frac{K_i}{2\pi\sinh r_i}$, and $-K_i\leq-\pi$ when $a\leq t\leq t_0$, then $r'_i(t)<0$ and hence $r_i(t)\leq r_i(a)=l$, which contradicts $r_i(t_0)>l$.

The arguments above show that $\overline{\{r(t)|t\in [0,+\infty)\}}$ is compactly contained in $\mathds{R}^N_{>0}$. We show there exists a sequence $\xi_n\uparrow +\infty$, such that $r(\xi_n)\rightarrow r^*$. In fact, $F(u(+\infty))$ exists. Furthermore, there exists a sequence $\xi_n\rightarrow +\infty$ such that
$$F(u(n+1))-F(u(n))=(F(u(t)))'|_{\xi_n}=\nabla F\cdot \frac{du}{dt}|_{\xi_n}=-\sum_i2\pi\sinh^2r_i\cdot \widetilde{R}_i^2 |_{\xi_n}\rightarrow 0.$$
Since $\{r(t)\}\subset\subset \mathds{R}^N_{>0}$, $r(\xi_n)$ have a convergent subsequence which is still denoted as $r(\xi_n)$ with $r(\xi_n)\rightarrow r^*$. Then we get $\widetilde{R}(r^*)=0$, and hence $K(r^*)=0$ and $R(r^*)=0$.
By Thurston's injective property of $K$-curvature map, $r^*$ is the unique zero curvature metric.

By calculating the differential of the right hand of (\ref{def-ricci flow2}) at $u^*$, we get
$$D|_{u^*}(-\widetilde{R})=-\frac{1}{2\pi}L\Sigma^{-1}$$
where $\Sigma=diag\{\sinh^2r_1, \cdots, \sinh^2r_N\}$. As $-L\Sigma^{-1}=-\Sigma^{\frac{1}{2}}\Sigma^{-\frac{1}{2}}L\Sigma^{-\frac{1}{2}}\Sigma^{-\frac{1}{2}}\sim \Sigma^{-\frac{1}{2}}L\Sigma^{-\frac{1}{2}}$, $D|_{u^*}(-\widetilde{R})$ has $N$ negative eigenvalues, which implies that $u^*$ is a local attractor of (\ref{def-ricci flow2}). Combining this with the fact that $u(\xi_n)$ converges to $u^*$, we get $u(t)\rightarrow u^*$ for any initial metric $u(0)$ with exponential convergence rate by ODE theory.\qed\\

\section{Convergence with nonpositive initial curvature metrics}\label{Convergence with nonpositive initial curvature metrics}
In the last subsection we proved combinatorial Ricci flow converges exponentially fast to zero curvature metric
if there exists a zero curvature metric.
It's interesting that we can get a convergence result without the assumption of the existence of zero curvature metric.
The main tool used in this subsection is a discrete version of maximum principle first developed by the authors in \cite{GX3}.

We first derive the evolution of modified curvature $\widetilde{R}_i$ along the combinatorial Ricci flow.

\begin{lemma}
Along the combinatorial Ricci flow (\ref{def-ricci flow}), the modified combinatorial curvature $R_i$ evolves according to
\begin{equation}\label{evolution of R}
\frac{d \widetilde{R}_i}{dt}=-\frac{1}{2\pi\sinh^2r_i}(L\widetilde{R})_i+2\cosh r_i\widetilde{R}_i^2.
\end{equation}
\end{lemma}
\proof
First note that $u_i=\ln \tanh\frac{r_i}{2}$, we have
$$\frac{\partial}{\partial u_i}=\sinh r_i\frac{\partial}{\partial r_i}.$$
Then we have
\begin{equation*}
\frac{\partial \widetilde{R}_i}{\partial u_j}
=\frac{1}{2\pi\sinh^2r_i}\frac{\partial K_i}{\partial u_j}-\frac{K_i}{2\pi\sinh^4r_i}\sinh r_j\frac{\partial}{\partial r_j}\sinh^2r_i
=\frac{1}{2\pi\sinh^2r_i}\frac{\partial K_i}{\partial u_j}-2\widetilde{R}_i\sinh r_j\coth r_i\delta_{ij},
\end{equation*}
which implies that
\begin{equation*}
\frac{d \widetilde{R}_i}{dt}
=\sum_j\frac{\partial \widetilde{R}_i}{\partial u_j}\frac{du_j}{dt}
=-\frac{1}{2\pi\sinh^2r_i}\sum_j\frac{\partial K_i}{\partial u_j}\widetilde{R}_j+2\cosh r_i\widetilde{R}_i^2.
\end{equation*}\qed

To get the estimation of curvature $R_i$, we can use the following discrete version of maximal principle introduced in \cite{GX3}.
\begin{proposition}\label{Maximum priciple}(Maximum Principle)
Let $f: V\times [0, T)\rightarrow \mathds{R}$ be a $C^1$ function such that
$$\frac{\partial f_i}{\partial t}\geq \Delta f_i+ \Phi_i(f_i), \ \ \forall (i, t)\in V\times [0, T) $$
where $\Delta f_i=\sum_{j\sim i}a_{ij}(f_j-f_i)$ with $a_{ij}>0$ and
$\Phi_i: \mathds{R}\rightarrow \mathds{R}$ is a local Lipschitz function.
Suppose there exists $C_1\in \mathds{R}$ such that $f_i(0)\geq C_1$ for all $i\in V$. Let $\varphi$ be the
solution to the associated ODE
\begin{equation*}
\begin{aligned}
\left\{
  \begin{array}{ll}
    \frac{d \varphi}{dt}=\Phi_i(\varphi) \\
    \varphi(0)=C_1,
  \end{array}
\right.
\end{aligned}
\end{equation*}
then
$$f_i(t)\geq \varphi(t)$$
for all $(i, t)\in V\times [0, T)$ such that $\varphi(t)$ exists.

Similarly, suppose $f: V\times [0, T)\rightarrow \mathds{R}$ be a $C^1$ function such that
$$\frac{\partial f_i}{\partial t}\leq \Delta f_i+ \Phi_i(f_i), \ \ \forall (i, t)\in V\times [0, T). $$
Suppose there exists $C_2\in \mathds{R}$ such that $f_i(0)\leq C_2$ for all $i\in V$. Let $\psi$ be the
solution to the associated ODE
\begin{equation*}
\begin{aligned}
\left\{
  \begin{array}{ll}
    \frac{d \psi}{dt}=\Phi_i(\psi) \\
    \psi(0)=C_2,
  \end{array}
\right.
\end{aligned}
\end{equation*}
then
$$f_i(t)\leq \psi(t) $$
for all $(i, t)\in V\times [0, T)$ such that $\psi(t)$ exists.
\end{proposition}

Using the Maximum principle, we have the following corollary.

\begin{corollary}\label{estimate of R by MP}
If $R_i(0)\leq 0$ for all $i\in V$, then $R_i(t)\leq 0$ for all $i\in V$ along the combinatorial Ricci flow (\ref{def-ricci flow}).
If $R_i(0)\geq 0$ for all $i\in V$, then $R_i(t)\geq 0$ for all $i\in V$ along the combinatorial Ricci flow (\ref{def-ricci flow}).
\end{corollary}
\proof
As $R_i$ and $\widetilde{R}_i$ always have the same sign, we just need to prove the corollary for $\widetilde{R}_i$.
As $L=A+L_B$ by Lemma \ref{Lemma-L-positive-definite}, we can rewrite the evolution equation for $\widetilde{R}$ as
$$\frac{d \widetilde{R}_i}{dt}=\Delta \widetilde{R}_i-\frac{A_i}{2\pi\sinh^2r_i}\widetilde{R}_i+2\cosh r_i\widetilde{R}_i^2,$$
where $\Delta \widetilde{R}_i=\sum_{j\sim i}a_{ij}(\widetilde{R}_j-\widetilde{R}_i)$ with $a_{ij}=\frac{B_{ij}}{2\pi\sinh^2r_i}>0$.
Applying the maximum principle gives the proof of the corollary. \qed

\begin{theorem}
Given $(X, T, \Phi,\mathds{H}^2)$. If the initial curvature are all nonpositive, i.e. $R_i(0)\leq0$ for each $i\in V$,
then there exists a zero curvature metric.
Furthermore, the solution to combinatorial Ricci flow (\ref{def-ricci flow})
exists for all time and converges to zero curvature metric exponentially fast.
\end{theorem}
\proof
Let $r(t)$ be the solution. The uniform upper bound of $r(t)$ follows the proof of Theorem  \ref{Thm-0curv-metric-exist-imply converg}.
For the lower bound, note that the negativity of the modified curvature $\widetilde{R}_i(t)$ is preserved along the flow.
Thus we have $\frac{dr_i}{dt}=-\widetilde{R}_i\sinh r_i\geq 0$, which implies that $r_i(t)\geq r_i(0)>0$. Hence we get $\{r(t)\}\subset\subset \mathds{R}^N_{>0}$.
Using Lemma \ref{Lemma-cpt-zcm-exists-converg} below, we get the proof.
\qed

\begin{lemma}\label{Lemma-cpt-zcm-exists-converg}
If the solution $r(t)$ to the combinatorial Ricci flow (\ref{def-ricci flow}) stays in a compact subset of $\mathds{R}^N_{>0}$, then there exists a unique zero curvature metric. Moreover, $r(t)$ converges exponentially fast to zero curvature metric.
\end{lemma}
\proof
The condition $\{r(t)\}\subset\subset \mathds{R}^N_{>0}$ implies the long time existence of the combinatorial Ricci flow. Consider the modified potential
\begin{equation}\label{def-G-ricci-potential}
G(u)\triangleq\int_{u_0}^u\sum_{i=1}^N K_idu_i,
\end{equation}
where $u_0$ is an arbitrary fixed point in $\mathds{R}^N_{<0}$. Then along the Ricci flow (\ref{def-ricci flow}), we have
\begin{equation}
\frac{dG(u(t))}{dt}=\nabla G\cdot \frac{du}{dt}=\sum_iK_i\cdot (-\widetilde{R}_i)=-\sum_i\frac{K_i^2}{2\pi\sinh r_i^2}\leq 0.
\end{equation}
So $G$ is decreasing along the flow (\ref{def-ricci flow}). Generally, $G$ is not bounded from below. However, $G(u(t))$ is bounded since $\overline{\{u(t)|t\in [0,+\infty)\}}$ is compact. Hence $G(u(+\infty))$ exists.
Then there exists a sequence $\xi_n\rightarrow +\infty$ such that
$$G(u(n+1))-G(u(n))=(G(u(t)))'|_{\xi_n}=\nabla G\cdot \frac{du}{dt}|_{\xi_n}=-\sum_i2\pi\widetilde{R}_i^2\cdot \sinh^2r_i|_{\xi_n}\rightarrow 0.$$
As $r_i(t)$ is bounded from above and away from zero, we have $\widetilde{R}_i(\xi_n)\rightarrow 0$.
Up to a subsequence, we can suppose that $r(\xi_n)\rightarrow r^*$ which has zero curvature.
We can finish the proof by showing that $D|_{u^*}(-\widetilde{R})=-\frac{1}{2\pi}L\Sigma^{-1}$ has $N$ negative eigenvalues and
hence $u^*$ is a local attractor of flow (\ref{def-ricci flow2}) or
by using the conclusion of Theorem \ref{Thm-0curv-metric-exist-imply converg} directly. \qed\\

\section{Combinatorial Ricci flow for general ``area element"}\label{Combinatorial Ricci flow for general area element}
For arbitrary $N$ positive functions $A_i: \mathds{R}^N_{>0}\rightarrow (0, +\infty)$, $r\mapsto A_i(r)$, we may take $A_i$ as the most general form of ``area element" at vertex $i$, where $1\leq i\leq N$. We shall study combinatorial Gaussian curvature and combinatorial Ricci flow for this generalized ``area element" $A_i$ in the following.

\begin{definition} ($A$-curvature) The modified combinatorial Gaussian curvature with respect to area element $A_i$ is defined as
\begin{equation}\label{def-A-curvature}
R_i=\frac{K_i}{A_i}.
\end{equation}
In the following we call it $A$-curvature for short.
\end{definition}

Obviously, a metric with zero $A$-curvature is exactly a metric with zero $K$-curvature. Thus by \cite{T1, CL1}, the zero $A$-curvature metric is unique if it exists. Inspired by the form (\ref{def-ricci flow2}) of combinatorial Ricci flow,
we introduce the following $A$-flow to study $A$-curvature.

\begin{definition} ($A$-flow) The hyperbolic combinatorial Ricci flow with respect to area element $A_i$ is defined as
\begin{equation}\label{def-A-curvature flow}
\frac{du_i}{dt}=-R_i,
\end{equation}
where $u_i=\ln \tanh \frac{r_1}{2}$. In the following we call it $A$-curvature flow or $A$-flow for short.
\end{definition}

\begin{theorem}
Given $(X, T, \Phi,\mathds{H}^2)$, $A$-flow (\ref{def-A-curvature flow}) converges if and only if there exists a zero curvature metric.
Furthermore, if the solution of (\ref{def-A-curvature flow}) converges, then it converges exponentially fast to the zero curvature metric.
\end{theorem}
\proof
For the ``only if" part, the proof is the same as that of Proposition \ref{prop-converg-imply-zccpm-exist}.
For the ``if" part, the methods is the same as in the proof of Theorem \ref{Thm-0curv-metric-exist-imply converg}.
We only list some differences here. In the proof of Theorem \ref{Thm-0curv-metric-exist-imply converg},
we first proved that all $r_i(t)$ are uniformly bounded below by a positive constant.
This procedure can be used here again to get a uniform lower bound of $r(t)$ as $A_i$ is always positive.

Suppose $T$ is the maximal existing time of $r(t)$ with $0<T\leq+\infty$. We claim that $r_i(t)$ is uniformly bounded above for all $i\in V$ and $t\in [0,T)$. If not, there exists $i\in V$ and $t_n\uparrow T$, such that $r_i(t_n)\uparrow +\infty$. Choose the same $l$ as in the proof of Theorem \ref{Thm-0curv-metric-exist-imply converg} and a $n_0$ such that $r_i(t_{n_0})>l$. Still denote $a=\inf\{t<t_{n_0}|r_i(t)>l\}$, then $r_i(a)=l$. Note that $A$-flow satisfies $r'_i(t)=-K_i \sinh r_i/A_i<0$, then $r_i(t)\leq r_i(a)=l$, which contradicts $r_i(t_{n_0})>l$. Thus $r_i(t)$ is uniformly bounded above for all $i$ and $t$, which implies that $T=+\infty$ and $\overline{\{r(t)|t\in [0,+\infty)\}}\subset\subset\mathds{R}^N_{>0}$. We can finish the proof by using Lemma \ref{Lemma-cpt-zcm-exists-converg}, which is still valid when we substitute $\sinh r_i^2$ by $A_i$. \qed\\

\textbf{Acknowledgements}\\[8pt]
The authors would like to thank Dr. Wenshuai Jiang for reading the paper carefully and 
giving some valuable suggestions to improve the writing of the paper. 
The first author would also like to give special thanks to Dr. Yurong Yu for her supports and encouragements during the work.
The research of the second author is partially supported by National Natural Science Foundation of China under
grant no. 11301402 and 11301399.
He would also like to thank Professor Guofang Wang for invitation to the Institute
of Mathematics of the University of Freiburg and for his encouragement and many
useful conversations during the work.

(Huabin Ge)
Department of Mathematics, Beijing Jiaotong University, Beijing 100044, P.R. China

E-mail: hbge@bjtu.edu.cn\\[2pt]

(Xu Xu) School of Mathematics and Statistics, Wuhan University, Wuhan 430072, PR China

E-mail: xuxu2@whu.edu.cn\\[2pt]

\end{document}